\documentclass[a4paper,10pt]{amsart}
\usepackage[arrow,matrix]{xy}
\usepackage{amsmath,amssymb,amscd,bbm,amsthm,mathrsfs,dsfont,enumerate}
\theoremstyle{plain} \textwidth=31pc \textheight=51pc

\usepackage{fancyhdr} 
\newtheorem{theorem}{Theorem}[section]
\newtheorem{lemma}[theorem]{Lemma}
\newtheorem{example}[theorem]{Example}
\newtheorem{proposition}[theorem]{Proposition}
\newtheorem{corollary}[theorem]{Corollary}
\theoremstyle{definition}
\newtheorem{definition}[theorem]{Definition}
\newtheorem{remark}[theorem]{Remark}
\numberwithin{equation}{section}

 \DeclareMathOperator{\Ext}{Ext}

\DeclareMathOperator{\e}{\epsilon}
\DeclareMathOperator{\GK }{GK dim}

\DeclareMathOperator{\ad}{ad}

\begin{document}

\title{Existence  of Hopf subalgebras of GK-dimension two} 
\author{Guangbin Zhuang}
\maketitle
\begin{abstract}
Let $H$ be a pointed Hopf algebra over an algebraically closed field of characteristic zero. If $H$ is a domain with finite Gelfand-Kirillov dimension greater than or equal to two, then $H$ contains a Hopf subalgebra of Gelfand-Kirillov dimension two.
\end{abstract}

\section{Introduction}
The study of Hopf algebras of low Gelfand-Kirillov dimension (or GK-dimension for short) provides profound insight into basic properties of general Hopf algebras. Noetherian prime regular Hopf algebras of GK-dimension one were studied in \cite{BZ, LWZ, G}.
Recently, in their paper \cite{GZ}, K. R. Goodearl and J. J. Zhang have classified all noetherian Hopf algebras $H$ over an algebraically closed field $k$ of characteristic zero, which are the integral domains of GK-dimension two and satisfy the condition $\Ext^1_H(k,k)\neq 0$.
After all the works listed above, it is natural to consider Hopf algebras of GK-dimension three. One way to understand a Hopf algebra with a higher GK-dimension is to find some of its Hopf subalgebras that are more familiar to us. As a small step toward this end, we will prove in this paper that for certain pointed Hopf algebras with GK-dimension $\ge 2$, there always exist a Hopf subalgebra of GK-dimension two.
We hope the result would offer some help in classifying the Hopf algebras of GK-dimension three.

The key step in achieving our main result is the analysis of the free pointed Hopf algebra $F(t)$ (see Example $\ref{original}$ for the definition).
In the classification result of \cite{GZ}, some of the Hopf algebras discussed there are the Hopf quotients of $F(t)$. When $t\ne 0$, $F(t)$ is a typical example of infinite-dimensional pointed non-cocommutative Hopf algebras. On the other hand, if $H$ is a pointed non-cocommutative Hopf algebra and $y\in H$ is a $(1, g)$-primitive element, then the Hopf subalgebra $K$ generated by $y$ and $g^{\pm 1 }$ is a Hopf quotient of $F(1)$. So it would be interesting to know more about the Hopf structure of $F(t)$.

In Section $3$, it is shown that $F(t)$, when $t\ne0$, is coradically graded in the sense of \cite[Lemma 2.2]{CM} with respect to the $y$-grading. In the same section, we also determine all skew-primitive elements of $F(t)$ when $t\ne 0$. It turns out that the coalgebra structure of $F(t)$, when $t\ne 0$, is quite different from that of $F(0)$. 

By understanding the Hopf algebra structure of $F(t)$, we obtain the following theorem on the existence of certain Hopf subalgebras of pointed Hopf algebras. Following the conventions in \cite{Mont}, for any Hopf algebra $H$, we denote the set of group-like elements in $H$ by $G(H)$, and the coradical filtration on $H$ by $\{H_n\}_{n\ge0}$. Also, given a group $G$, we will use the standard notation $k[G]$ to denote the group algebra of $G$ over $k$.

\begin{theorem}[Theorem $\ref{GK2sub}$]
Let $H$ be a pointed Hopf algebra over an algebraically closed field $k$ of characteristic $0$. If $2\le \GK H<\infty$ and $H$ is a domain, then $H$ contains a Hopf subalgebra of GK-dimension $2$. To be precise, 
\begin{enumerate}
\item[\textup{(1)}]  If $\GK k[G]=0$, where $G=G(H)$, then $H$ contains a Hopf subalgebra isomorphic to $U(\mathbf{h})$, where $\mathbf{h}$ is a $2$-dimensional Lie algebra.
\item[\textup{(2)}]  If $\GK k[G]\ge 2$, then $H$ contains a Hopf subalgebra isomorphic to $k\mathbb{Z}^2$.
\item[\textup{(3)}]  If $H_1\supsetneq H_0=k[G]$ and $\GK k[G]\ge 1$, then $H$ contains a Hopf subalgebra isomorphic to either $A(b, \xi)$ or $C(b+1)$ (see Example $\ref{typeA}$ and $\ref{typeC}$) for some integer $b$ and $\xi \in k^{\times }$.
\end{enumerate}
\end{theorem}
In fact, the characteristic zero condition is only used in part $(1)$. Also, part $(1)$ and $(2)$ of Theorem $1.1$ do not require $H$ to be a domain. Moreover, as suggested by known examples, it is reasonable to conjecture that Theorem $1.1$ still holds if $H$ is weakened to a prime algebra, except that in part $(3)$ we might have to include some other classes of Hopf algebras, besides $A(b, \xi)$ and $C(b+1)$, as possible Hopf subalgebras.

The essential ingredient for proving the main result is the following theorem.
\begin{theorem}[Corollary \ref{existence2}]\label{intro}
Let $K$ be a pointed Hopf algebra with finite GK-dimension. If $y\in K_1\backslash K_0$ is a $(1, g^b)$-primitive element for some group-like element $g\ne1$ and integer $b$, then there are $b_0, \cdots, b_m, \beta\in k$ and $\xi \in k^{\times }$ such that $f=b_0y+b_1gyg^{-1}+\cdots+b_{m}g^{m}yg^{-m}\in  K_1\backslash K_0$ and $gf=\xi fg+\beta g(g^b-1)$.
\end{theorem}

In fact, Theorem $\ref{intro}$ can be viewed as an infinite-dimensional analogue of the following theorem of D. \c{S}tefan.

\begin{theorem}\cite[Theorem 2]{S}\label{Stefan}
Let $H$ be a finite-dimensional pointed Hopf algebra. If $H$ is not cosemisimple, then there exist two natural numbers $m$ and $n$, with $m\ne 1$ and $m$ dividing $n$, an $m$-th primitive root $\omega$ of unity and two elements $g, x\in H$ such that
\begin{enumerate}
\item[\textup{(1)}] $gx=\omega xg$;
\item[\textup{(2)}] $g$ is a group-like element of order $n$;
\item[\textup{(3)}] $x\in P_{1, g}$ and $x^m$ is either $0$ or $g^m-1$.
\end{enumerate}
\end{theorem}
With the notation in Theorem $\ref{Stefan}$, it is easy to check that the Hopf subalgebra generated by $g$ and $x$ is a Hopf quotient of $A(1, \omega)$.

\section{Preliminaries}
In this section, we give some basic definitions and properties related to coalgebras, which will be intensively used in the following sections.
Throughout this paper, $k$ always denotes a base field which is algebraically closed of characteristic $0$. All algebras, coalgebras and tensor products are taken over $k$. Given a group $G$, we will use the standard notation $k[G]$ to denote the group algebra of $G$ over $k$.

Suppose that $C$ is a coalgebra. Let $G(C)$ be the set of group-like elements in $C$. An element $y\in C$ is called \textbf{$(g, h)$-primitive} if $g, h\in G(C)$ and $\Delta(y)=y\otimes g+h\otimes y$. We use $P_{g, h}(C)$ to denote the set of $(g, h)$-primitive elements in $C$ \cite[p. 67]{Mont}. If $C=H$ is a Hopf algebra, we simply write $P(H)$ for $P_{1, 1}(H)$, the set of primitive elements in $H$. 

The \textbf{coradical}  $C_0$ of $C$ is defined to be the sum of all simple subcoalgebras of $C$. The coalgebra $C$ is called \textbf{cosemisimple} if $C=C_0$, \textbf{pointed} if $C_0=kG(C)$, and \textbf{connected} if $C_0$ is one-dimensional.

We now recall the definitions of graded coalgebras and coalgebra filtrations, which will play a central role in the following sections.
\begin{definition}
Let $C$ be a coalgebra. If 
\begin{enumerate}
\item$C=\bigoplus\limits_{i=0}^{\infty }C(i)$ as $k$-vector spaces,
\item$\Delta(C(n))\subset\bigoplus\limits_{i=0}^{n}C(i)\otimes C(n-i) $ for any $n\ge0$,
\item $\e(C(n))=0$ for any $n\ge 1$,
\end{enumerate}
then $C=\bigoplus\limits_{i=0}^{\infty }C(i)$ is called a \textbf{graded coalgebra}.
\end{definition}

\begin{definition}
Let $C$ be a coalgebra and $\{A_n\}_{n\ge 0}$ a family of subcoalgebras of $C$. Then $\{A_n\}_{n\ge 0}$ is called a \textbf{coalgebra filtration} if
\begin{enumerate}
\item $A_n\subset A_{n+1}$ and $C=\cup_{n\ge 0} A_n$,
\item$\Delta(A_n)\subset\sum\limits_{i=0}^{n}A_i\otimes A_{n-i} $ for any $n\ge0$,
\end{enumerate}
\end{definition}

Starting from $C_0$, one can build up a canonical coalgebra filtration $\{C_n\}$ of $C$, which is defined inductively by 
$$C_n=\Delta^{-1}(C\otimes C_{n-1}+ C_0\otimes C).$$
This filtration is called the \textbf{coradical filtration} of $C$. 

The following lemma is actually in the proof of \cite[Theorem 5.3.1]{Mont}.
\begin{lemma}\label{ghprimitive}
Let $C$ be a coalgebra with coradical filtration $\{C_n \}$ and let $I$ be a non-zero coideal of $C$. Then $I\cap C_1\ne \{0\}$. Moreover, if $C$ is pointed, then $I$ contains a non-zero $(g, h)$-primitive element for some $g, h\in G(C)$.
\end{lemma}
\begin{proof}
Consider the coalgebra map $\pi : C\rightarrow C/I$. Then the kernel of $\pi \mid_{C_1}$ is $I\cap C_1$. Suppose $I\cap C_1=\{ 0\}$, then $\pi |_{C_1}$ is injective. By \cite[Theorem 5.3.1]{Mont}, $\pi$ is injective. But this implies $I=\{ 0\}$, which contradicts the choice  of $I$.

Suppose further that $C$ is pointed. The kernel of $\pi |_{P_{g, h}(C)}$ is $I\cap P_{g, h}(C)$. Suppose $I\cap P_{g, h}(C)=\{ 0\}$ for any $g, h\in G(C)$, then $\pi |_{P_{g, h}(C)}$ is injective for any $g, h\in G(C)$. By \cite[Corollary 5.4.7]{Mont}, $\pi$ is injective, which is a contradiction.
\end{proof}

Given a graded coalgebra $C=\bigoplus\limits_{i=0}^{\infty }C(i)$, one can define a coalgebra filtration $\{A_n\}_{n\ge 0}$ by setting $A_n=C(0)\oplus\cdots \oplus C(n)$. The Hopf algebra $F(t)$ (see Example \ref{original} for the definition), which we will study in the following sections, has a natural graded coalgebra structure that is compatible with the coradical filtration in the following sense.

\begin{definition}
Let $C=\bigoplus\limits_{i=0}^{\infty }C(i)$ be a graded coalgebra with coradical filtration $\{C_n\}$. Then $C$ is called \textbf{coradically graded} if $C_0=C(0)$ and $C_1=C(0)\bigoplus C(1)$. A graded Hopf algebra $H$ is coradically graded if it is coradically graded as a graded coalgebra.
\end{definition}

The terminology is justified by the following lemma.
\begin{lemma}\cite[Lemma 2.2]{CM}
If $C$ is coradically graded, then 
$$C_m=\bigoplus\limits_{i\le m }C(i).$$
\end{lemma}

\section{Free Pointed Hopf algebra $F(t)$}
In this section we introduce the free pointed Hopf algebra $F(t)$ and study some of its basic properties.
\begin{example}\label{original}
Let $H=k\langle x^{\pm 1}, y \rangle$ be the algebra with relations $xx^{-1}=x^{-1}x=1$, then $H$ becomes a Hopf algebra via
$$\Delta(x)=x\otimes x, \,\,\,\Delta (y)=y\otimes1+x^t\otimes y,$$
$$S(x)=x^{-1}, \,\,\,S(y)=-x^{-t}y,$$
$$\epsilon(x)=1, \,\,\, \epsilon(y)=0,$$
where $t$ is any given integer. We denote this Hopf algebra by $F(t)$.
\end{example}

\begin{proposition}\label{basis}
The Hopf algebra $F(t)$ has a $k$-linear basis 
$$\{x^{i_1}y x^{i_2}y\cdots yx^{i_{n+1}} \mid\,\,(i_1,\cdots, i_{n+1})\in \mathbb{Z}^{n+1}, \,\,n\in \mathbb{N}\}.$$
\end{proposition}
\begin{proof}
This is standard from the diamond lemma.
\end{proof}

The Hopf algebra $F(t)$ has several gradings, which turn out to be crucial in determining its coradical filtration. 
\begin{proposition}\label{ygrading}
There is a natural $y$-grading on $F(t)$ under which $F(t)$ becomes a graded coalgebra.
\end{proposition}
\begin{proof}
The natural $y$-grading is given by setting $\deg_yy=1$ and $\deg_yx^{\pm 1 }=0$. It is easy to check that $F(t)$ becomes a graded algebra by this setting, since the relations all preserve the grading. This $y$-grading on $F(t)$ induces a natural $y$-grading on $F(t)\otimes F(t)$ by defining $\deg_y(f_1\otimes f_2)= \deg_y f_1+\deg_y f_2$, where $f_1, f_2\in F(t)$. Also notice that this grading on $F(t)\otimes F(t)$ is an algebra grading.

By definition of a graded coalgebra, we only have to show that if $h\in F(t)$ is homogeneous of $y$-degree $n$, then so is $\Delta(h)$. By Proposition \ref{basis}, we may assume that $h$ is of the form $x^{i_1}y x^{i_2}y\cdots yx^{i_{n+1}} $, where $(i_1,\cdots, i_{n+1})\in \mathbb{Z}^{n+1}$.
Now
$$\Delta(h)=\Delta(x^{i_1})\Delta(y) \cdots \Delta(y)\Delta(x^{i_{n+1}}). $$
Notice that $\Delta(x^{i_j})$ and $\Delta(y)$ are homogeneous of $y$-degrees $0$ and $1$, respectively. The result then follows from the fact that $F(t)\otimes F(t)$ is a graded algebra with respect to the $y$-grading. This completes the proof.
\end{proof}

By analyzing the grading, we are actually able to determine the coradical of $F(t)$. Notice that by construction, the degree zero part of $F(t)$ with respect to the $y$-grading is just $k[x^{\pm 1}]$, the subalgebra of $F(t)$ generated by $x^{\pm1}$.

\begin{corollary}\label{grouplike}
Let $H=F(t)$. Then $H_0=k[x^{\pm1 }]$. As a consequence, $F(t)$ is a pointed Hopf algebra.
\end{corollary}
\begin{proof}
By the previous proposition, $H$ is a graded coalgebra with respect to the $y$-grading. The coalgebra grading induces a coalgebra filtration $\{A_n\}_{n\ge 0}$ on $H$ with $A_0=k[x^{\pm 1}]$. By \cite[Lemma 5.3.4]{Mont}, $A_0\supset H_0$. On the other hand, $A_0$ is spanned by group-like elements and thus $k[x^{\pm1 }]=A_0=H_0$.
\end{proof}

There is also a natural $x$-grading on $F(t)$ given by setting $\deg_xx^{\pm 1 }=\pm1$ and $\deg_x{y}=0$.

\begin{remark} Note that $F(0)$ is a cocommutative Hopf algebra. In this case, the Hopf structure is well understood. In fact, by \cite[Corollary 5.6.4, Theorem 5.6.5]{Mont}, $F(0)\cong U(\mathbf{g})\#k[x^{\pm1 }]$, where $\mathbf{g}$ is the set of primitive elements in $F(0)$. Actually, $\mathbf{g}$ is the free Lie algebra generated by $\{ x^nyx^{-n}\mid n \in \mathbb{Z}\}$. It is easy to see that $F(0)$ is coradically graded by setting $\deg x^{\pm 1}=0$ and $\deg f=1$ for any $f\in \mathbf{g}$. Moreover, it is not hard to prove that every non-zero Hopf ideal $I$ of $F(0)$ is generated by elements in $\mathbf{g}$ and possibly an element of the form $x^n-1$ for some positive integer $n$. So in next section, we will focus on the case $t\ne 0$.
\end{remark}


\section{The coradical filtration of $F(t)$ when $t\ne 0$}
In the section, we study the coradical filtration of the free Hopf algebra $F(t)$ described in Example \ref{original} with $t\ne0$. We can further assume that $t>0$, since it is obvious that $F(t)\cong F(-t)$ as Hopf algebras for any integer $t$. We will show that $H$ is a coradically graded Hopf algebra with respect to the $y$-grading. This is quite different from the $t=0$ case. In fact, if $t=0$, then $P(F(0))$ is the free Lie algebra generated by $\{ x^nyx^{-n}\mid n \in \mathbb{Z}\}$, which implies that there are primitive elements of any $y$-degree in $F(0)$.

Throughout this section, let $H$ be the Hopf algebra $F(t)$ with $t\ne0$.
The following proposition is crucial in proving that $H$ is coradically graded with respect to the $y$-grading. For the rest of this section, we use $H(n)$ to denote the $k$-linear subspace spanned by homogeneous elements of $y$-degree $n$.

\begin{lemma}\label{nocancel}
Any skew-primitive element in $H$ has $y$-degree $\le 1$.
\end{lemma}

\begin{proof}
Suppose to the contrary that $f$ is an $(x^a, x^c)$-primitive element of $y$-degree $n\ge2$. Replacing $f$ with $x^{-a}f$, we can assume that $f$ is $(1, x^b)$-primitive for some integer $b$, that is to say,
\begin{equation}\label{comul1}
\Delta(f)=f\otimes 1+ x^b\otimes f. 
\end{equation}
It is also clear from Proposition $\ref{ygrading}$ that  we can further assume that $f$ is homogeneous of $y$-degree $n$.

For any $\alpha=(i_1,\cdots, i_{n+1})\in \mathbb{Z}^{n+1}$, denote the monomial $x^{i_1}y x^{i_2}y\cdots yx^{i_{n+1}} $ by $M_{\alpha }$ and $\sum\limits_{s=1}^{n} i_{s}$ by $T(\alpha)$. Then we can write
$$f=\sum \limits_{\alpha\in \mathcal{I} }\mu_{\alpha }M_{\alpha },$$
where $\mathcal{I}$ is a non-empty finite subset of $\mathbb{Z}^{n+1}$ and $\mu_{\alpha }\in k^{\times }$ for any $\alpha\in \mathcal{I}$. A direct computation shows that 
$$\Delta(M_{\alpha })=M_{\alpha }\otimes x^{T(\alpha)}+x^{T(\alpha)+nt}\otimes M_{\alpha }+w_{\alpha },  $$
where $w_{\alpha }\in \sum\limits_{i, j< n \atop i+j=n }H(i)\otimes H(j)$. Hence 
\begin{equation}\label{comul2}
\Delta(f)=\sum\limits_{\alpha\in \mathcal{I} }\mu_{\alpha } (M_{\alpha }\otimes x^{T(\alpha)}+ x^{T(\alpha)+nt}\otimes M_{\alpha })+ w,
\end{equation}
where $w\in \sum\limits_{i, j< n \atop i+j=n }H(i)\otimes H(j)$. By comparing $(\ref{comul1})$ and $(\ref{comul2})$ we see that $T(\alpha )=0$ for any $\alpha\in \mathcal{I}$, and $b=nt$.

Define two bijective maps $\sigma$ and $\tau$ from $\mathbb{Z}^{n+1}$ to $\mathbb{Z}^{n+1}$ by
$$\sigma(i_1,\cdots, i_{n+1})=(-i_{n+1}-t,-i_{n}-t, \cdots, -i_2-t, -i_1),$$
$$\tau(i_1,\cdots, i_{n+1})=(i_1-b,i_2, \cdots, i_{n+1}).$$
Since the antipode $S$ is an anti-algebra map, it is easy to check that $S(f)=\sum \limits_{\alpha\in \mathcal{I} }(-1)^n\mu_{\alpha }M_{\sigma \alpha }$. On the other hand, since $S$ is the inverse to $Id_{H}$ under convolution, we have $S(f)=-x^{-b}f$, which implies that $S(f)=-\sum \limits_{\alpha\in \mathcal{I} }\mu_{\alpha }M_{\tau \alpha }$. It then follows that 
\begin{equation}\label{key}
\sum \limits_{\alpha\in \mathcal{I} }((-1)^n\mu_{\alpha }M_{\sigma \alpha }+M_{\tau \alpha })=0.
\end{equation}
Choose $\beta=(\ell_1,\cdots, \ell_{n+1})\in \mathcal{I}$. From $(\ref{key})$ we see that there is $\beta'\in \mathcal{I}$ such that $\tau\beta'=\sigma\beta$, i.e. $(\tau^{-1}\sigma)\beta\in \mathcal{I}$. By induction, the set $\mathcal{A}:=\{(\tau^{-1}\sigma)^i\beta\,\mid\, i\ge0 \}$ is a subset of $\mathcal{I}$. But it is easy to show that
$$(\tau^{-1}\sigma)^{2s}\beta=(\ell_{1}+s(n-1)t, \ell_2, \cdots, \ell_n, \ell_{n+1}-s(n-1)t)$$
for any $s\ge0$. Since $(n-1)t>0$, by comparing the first entry, it is clear that $(\tau^{-1}\sigma)^{2s}\beta\ne (\tau^{-1}\sigma)^{2s'}\beta$ whenever $s\ne s'$. Hence the set $\mathcal{A}$ is infinite, which is a contradiction. This completes the proof.
\end{proof}


\begin{proposition}\label{cofiltration}
$H$ is a coradically graded Hopf algebra with respect to the $y$-grading.
\end{proposition}
\begin{proof}
By Corollary $\ref{grouplike}$, $H_0=H(0)=k[x^{\pm 1}]$. Also it is easy to check that $H_1\supset H(0)\oplus H(1)$. By \cite[Theorem 5.4.1]{Mont}, $H_1$ is spanned by group-like and skew-primitive elements. By Lemma $\ref{nocancel}$, any skew-primitive element in $H$ has $y$-degree $\le 1$. Hence $H_1\subset H(0)\oplus H(1)$, which completes the proof.
\end{proof}

With this proposition, we can give a full description of all skew-primitive elements in $H$. For any $h\in H$, we define $\deg h$ to be the pair $(\deg_x h, \deg_y h)$. Notice that this $(\mathbb{Z}\times\mathbb{N})$-grading on $H$ induces a $(\mathbb{Z}\times\mathbb{N}\times\mathbb{Z}\times \mathbb{N})$-grading on $H\otimes H$.
\begin{proposition}\label{primitiveofH}
Let $h$ be a skew-primitive element in $H$. Then $\deg_y h \le 1$. If $\deg_y h=0$, i.e. $h$ is a polynomial in $x^{\pm 1 }$, then $h$ is of the form $\lambda x^a(x^m-1)$, where $ \lambda \in k^{\times }$, $a\in \mathbb{Z}$ and $m$ is a positive integer. If $\deg_y h=1$, then $h$ is of the form $f+\lambda x^a(x^t-1)$, where $\lambda\in k$ and $f$ is homogeneous with $\deg f= (a, 1)$.
\end{proposition}
\begin{proof}
By the definition of coradical filtration, $h\in H_1$. By Proposition $\ref{cofiltration}$, $H_1=H(0)\oplus H(1)$. So $\deg_y h \le 1$. 

If $\deg_y h=0$, then $h$ is a polynomial in $x^{\pm1 }$. Let $h=\sum \alpha_i x^i$ and suppose that $h$ is $(x^a, x^b)$-primitive for some $a, b\in \mathbb{Z}$. Then 
\begin{equation}\label{case1}
\Delta(h)=\Delta(\sum \alpha_i x^i)=\sum \alpha_i x^i\otimes x^i.
\end{equation}
On the other hand, since $h$ is $(x^a, x^b)$-primitive, it yields
\begin{equation}\label{case2}
\Delta(h)=h\otimes x^a+ x^b\otimes h.
\end{equation}
By comparing $(\ref{case1})$ and $(\ref{case2})$, we see that $\alpha_i =0$ if $i\neq a$ and $i\ne b$. Therefore, $h=\alpha_a x^a+ \alpha_bx^b$ and 
\begin{align}
\Delta(h)=&\Delta(\alpha_a x^a+ \alpha_bx^b)\\
=&(\alpha_a x^a+ \alpha_bx^b)\otimes x^a+x^b\otimes (\alpha_a x^a+ \alpha_bx^b)-\alpha_bx^b\otimes x^a-\alpha_ax^b\otimes x^a\nonumber \\
=&h\otimes x^a+ x^b\otimes h-(\alpha_a+\alpha_b)x^b\otimes x^a.\nonumber
\end{align}
It then follows that $\alpha_a=-\alpha_b\ne 0$. Hence $h$ has the claimed form.

If $\deg_y h=1$, we can write $h(x,y)=f_1+\cdots+f_s+ g$, where $f_i$ is homogeneous with $\deg f_i=(a_i, 1)$, and $g$ is a polynomial in $x^{\pm 1 }$. First, we need to show $s=1$. If not, we can rearrange $f_i$ so that $a_1>a_2>\cdots>a_s$. Suppose $h$ is an $(x^m, x^n)$-primitive element, that is, 
\begin{align}\label{compare1}
\Delta(h)=&h\otimes x^m+x^n\otimes h \\
   =&f_1\otimes x^m+x^n\otimes f_1+\cdots + f_s\otimes x^m+x^n\otimes f_s+\cdots\nonumber
\end{align}
On the other hand, it is easy to check that $f_i$ is an $(x^{a_i}, x^{a_i+t})$-primitive element, so
\begin{align}\label{compare2}
\Delta(h)=&f_1\otimes x^{a_1}+x^{a_1+t}\otimes f_1+\cdots + f_s\otimes x^{a_s}+x^{a_s+t}\otimes f_s+\cdots
\end{align}
Since $s>1$, there is some $j$ such that $a_j\ne m$. From $(\ref{compare2})$, we see that $ f_j\otimes x^{a_j}$ is the non-zero component of $\Delta(h)$ in degree $(a_j, 1;a_j, 0)$. But $(\ref{compare1})$ implies that $\Delta(h)$ has no component in degree $(a_j, 1;a_j, 0)$. This is a contradiction.

Now we can write $h=f+g$, where $f$ is homogeneous with $\deg f= (a, 1)$. If $g=0$, we are done. So assume $g\ne 0$ and it suffices to show that $g$ is of the form $\lambda x^a(x^t-1)$ for some $\lambda \in k^{\times }$. Note that $f$ is a $(x^a, x^{a+t})$-primitive element, so similar to the argument above, $h$ must be a $(x^a, x^{a+t})$-primitive element. Hence $g$ is $(x^a, x^{a+t})$-primitive. Then it follows from the previous case where $\deg_y h=0$ that $g$ is of the form $\lambda x^a(x^t-1)$, as desired.
\end{proof}

\begin{corollary}\label{existence}
Suppose that $I$ is a non-zero Hopf ideal of $H$. Then $I$ contains an element $f$ of the form $x^m-1$ for some integer $m$ or $a_0y + a_1xyx^{-1}+\cdots+a_nx^nyx^{-n}+\lambda(x^t-1)$, where $a_0, \cdots , a_n, \lambda\in k$ and $a_0, a_n\ne 0$.
\end{corollary}
\begin{proof}
Just combine Lemma $\ref{ghprimitive}$ and Proposition $\ref{primitiveofH}$.
\end{proof}

\begin{corollary}\label{existence2}
Let $K$ be a pointed Hopf algebra with finite GK-dimension. If $y\in K_1\backslash K_0$ is a $(1, g^b)$-primitive element for some group-like element $g\ne1$ and integer $b$, then there are $b_0, \cdots, b_m, \beta\in k$ and $\xi \in k^{\times }$ such that $f=b_0y+b_1gyg^{-1}+\cdots+b_{m}g^{m}yg^{-m}\in  K_1\backslash K_0$ and $gf=\xi fg+\beta g(g^b-1)$.
\end{corollary}
\begin{proof} 
First we claim that there are $a_0, \cdots , a_n, \lambda\in k$ with $a_0\ne 0, n\ge 1$ such that $a_0y + a_1gyg^{-1}+\cdots+a_ng^nyg^{-n}+\lambda(g^b-1)=0$. 
If $y$ is primitive, i.e. $g^b=1$, then the $k$-linear space $V$ spanned by $\{ g^iyg^{-i}\mid i\in\mathbb{Z}\}$ is a subspace of $\mathbf{g}$, the $k$-linear space spanned by all primitive elements in $K$. By the universal property of $U(\mathbf{g})$, the inclusion $\mathbf{g}\hookrightarrow K$ extends to an algebra map $\tau: U(\mathbf{g})\rightarrow K$. A direct calculation by using the PBW basis of $U(\mathbf{g})$ shows that $\tau$ is also a coalgebra map. By  \cite[Theorem 5.3.1, Proposition 5.5.3]{Mont}, $\tau$ is injective. It then follows from \cite[Lemma 3.1]{KL} that $\GK U(\mathbf{g})\le \GK K< \infty$. Now we have $\dim_k V\le \dim_k\mathbf{g}=\GK U(\mathbf{g})<\infty$. Hence $a_0y + a_1gyg^{-1}+\cdots+a_ng^nyg^{-n}=0$, for some $a_0, \cdots , a_n \in k$ with $a_0\ne 0$ and $n>0$. If $g$ is torsion with the order $n\ge2$, then $y-g^nyg^{-n}=0$. If $g$ is torsion-free and $b\ne0$, the Hopf subalgebra generated by $y$ and $g^{\pm 1 }$ is a Hopf quotient of $F(b)$ where $b\neq 0$. Now the statement follows from Corollary $\ref{existence}$. 

Choose $n$ to be minimal. Since $y\in K_1\backslash K_0$, we have $n\ge1$. Let $\xi$ be a root of the polynomial $a_0+ a_1s+\cdots+a_ns^n$, and set $b_0+b_1s+\cdots+b_{n-1}s^{n-1}=(a_0+ a_1s+\cdots+a_ns^n)/(s-\xi)$. Let $f=b_0y+b_1gyg^{-1}+\cdots+b_{n-1}g^{n-1}yg^{-(n-1)}$.  If $f\in K_0$, then $f\in P_{1, g^b}(K)\cap K_0=k(g^b-1)$. The equality $P_{1, g^b}(K)\cap K_0=k(g^b-1)$ can be verified by straight calculation since $K_0$ is spanned by group-like elements. Hence $f=\gamma (g^b-1)$ for some $\gamma \in k$, which contradicts the minimality of $n$. As a consequence, $f\in K_1\backslash K_0$. By the choice of $f$, we have $gfg^{-1}-\xi f+ \beta(g^b-1)=0$ for some $\beta \in k$. This completes the proof.
\end{proof}


\section{Main Theorem}
This section is devoted to the proof of the main theorem (Theorem \ref{GK2sub}). We also explain why the theorem fails if some of the conditions are weakened.
To begin with, we list two classes of Hopf algebras of GK-dimension $2$, which are defined in \cite[Construction 1.1 and 1.4]{GZ}. It turns out that these two classes of Hopf algebras appear naturally as Hopf subalgebras of a large class of Hopf algebras (see Lemma \ref{case3}).

\begin{example}\label{typeA}
Let $A=k\langle g^{\pm1 }, y \mid gy=\xi yg\rangle$, where $\xi \in k^{\times }$. Then $A$ has a Hopf algebra structure under which $x$ is group-like and $y$ is $(1, g^b)$-primitive. We denote it by $A(b, \xi)$. 
\end{example}

\begin{example}\label{typeC}
Let $C=k\langle g^{\pm1 }, y \mid gy=yg+g^{b+1}-g\rangle$. Then $C$ also has a Hopf algebra structure under which $g$ is group-like and $y$ is $(1, g^b)$-primitive. We denote it by $C(b+1)$.
\end{example}

\begin{remark}
The definition of $C(b+1)$ in Example \cite[Construction 1.4]{GZ} is slightly different from Example \ref{typeC}. However, it is not hard to show that they are isomorphic as Hopf algebras. 
In \cite[Construction 1.1 and 1.4]{GZ}, it is also shown that $A(m, r)\cong A(n, q)$ if and only if $(m, r)=(n, q)$ or $(m, r)=(-n, q^{-1})$, and that $C(m)\cong C(n)$ if and only if $m=n$. Moreover, both classes are domains of GK-dimension $2$.
\end{remark}

\begin{lemma}\label{AandC}
Let $H$ be a noetherian prime Hopf algebra of GK-dimension $2$. Suppose that $\pi : H\rightarrow K$ is a Hopf algebra projection to a domain $K$. If $K$ is neither cosemisimple nor connected, then $\pi$ is an isomorphism.
\end{lemma}
\begin{proof}
If $\pi$ is not an isomorphism, then $I:=\ker\pi \ne0$. By Goldie's Theorem, the ideal $I$ contains a regular element.  By \cite[Proposition 3.15]{KL}, $\GK K\le \GK H-1=1$. Since there is no algebra with GK-dimension strictly between $0$ and $1$, $\GK K=0$ or $1$. By assumption, $K$ is not finite-dimensional, so $\GK K=1$. By \cite[Lemma 4.5]{GZ}, $K$ is commutative. Hence $K$ is a commutative Hopf domain with $\GK K=1$. But by \cite[Proposition 2.1]{GZ}, $K$ is either a group algebra or an enveloping algebra. Hence $K$ is either cosemisimple or connected, which contradicts the assumption.
\end{proof}

The next three lemmas are known. However, I am not able to locate them in literature and thus the proofs are included.
\begin{lemma}\label{addition}
Let $G$ be a finitely generated group which acts on a finitely generated $k$-algebra $A$. Suppose that there is a finite-dimensional generating $k$-subspace $W$ of $A$ such that $W$ contains $1$ and is stable under the action of $G$. Then 
$$\GK A *G= \GK A + \GK k[G],$$
where $A*G$ is the skew group algebra.
\end{lemma}
\begin{proof} 
It suffices to assume that both $\GK A$ and $\GK k[G]$ are finite. By \cite[Theorem 11.1]{KL}, $G$ has a nilpotent normal subgroup $N$ with finite index. By Schreier's lemma, $N$ is finitely generated. By \cite[Proposition 5.5]{KL}, $\GK A*G=\GK A*N$ and $\GK k[G]=\GK k[N]$. Hence by replacing $G$ with $N$, we can assume $G$ is finitely generated nilpotent.

Suppose that the group $G$ is generated by the set $S=\{g_1, g_2, \cdots, g_m\}$. Let $g_S(n)$ be the number of distinct group elements that can be expressed as words of length $\le n$ in the specified generators and their inverses. Without loss of generality, we can assume that $S$ is closed under inversion. By the definition of GK-dimension and that of group algebras,  $\GK k[G]= \limsup\limits_{n\rightarrow \infty }\log_n g_S(n)$. Moreover, if $\GK k[G]< \infty$, then $\GK k[G]=\lim\limits_{n\rightarrow \infty }\log_ng_S(n)$ (see \cite[Theorem 11.14]{KL}).
Let $V=W+Wg_1+\cdots+ Wg_m$. It is clear that $V$ is finite-dimensional and generates $A *G$.

First, we prove that $\GK  A *G \ge\GK A + \GK k[G]$. Notice that 
$$\displaystyle \sum _i W^nf_i\subset (W+Wg_1+\cdots+ Wg_m)^{2n}=V^{2n},$$
where $f_i$ runs through all words of length $\le n$. So $g_S(n)\cdot \dim W^n\le \dim V^{2n}$. It then follows that 
\begin {align*}\displaystyle 
\GK A + \GK k[G] &=\limsup_{n\rightarrow \infty }\log_n W^n +\lim_{n\rightarrow \infty }\log_ng_S(n)\\
&=\limsup_{n\rightarrow \infty }\log_n(g_S(n)\cdot\dim W^n)\\
&\le\limsup_{n\rightarrow \infty }\dim V^{2n}\le \GK A *G.
\end{align*}
For the other direction, notice that 
$$\displaystyle (W+Wg_1+\cdots+ Wg_m)^{n}\subset  \sum _i W^nf_i,$$
where $f_i$ runs through all words of length $\le n$. Here we use the fact that $g_jW=Wg_j$ for any $j$. As a consequence, $ \dim V^{n}\le g_S(n)\cdot \dim W^n$. Therefore, 
\begin{align*}
\GK A *G=\limsup_{n\rightarrow \infty }\log_n\dim V^{n}&\le \limsup_{n\rightarrow \infty }\log_n(g_S(n)\cdot\dim W^n)\\
&\le \limsup_{n\rightarrow \infty }\log_n W^n +\lim_{n\rightarrow \infty }\log_ng_S(n)\\
&=\GK A + \GK k[G]. 
\end{align*}
This completes the proof.
\end{proof}

\begin{lemma}\label{groupalgebra}
Let $G$ be a group and $N$ a finite normal subgroup of $G$. Then 
$$\GK k[G]=\GK k[G/N].$$
\end{lemma}
\begin{proof}
Notice that there is a surjective algebra map $k[G]\rightarrow k[G/N]$. By \cite[Lemma 3.1]{KL}, $\GK k[G]\ge\GK k[G/N]$. Now we prove the other direction. Denote $G/N$ by $\overline{G}$ and $k[N]$ by $A$. By \cite[Example 7.1.6]{Mont}, 
$$k[G]\cong A\#_{\sigma}k[\overline{G}],$$
where $\sigma: k[\overline{G}]\otimes k[\overline{G}]\rightarrow A$ is a cocycle. The multiplication on $A\#_{\sigma}k[\overline{G}]$ is given by 
\begin{equation}\label{product}
(a\#g)(b\#h)=a(g.b)\sigma(g, h)\#gh,
\end{equation}
where $a, b\in N$ and $g, h\in \overline{G}$. Here $g.b$ is the image of $g\otimes b$ under a certain $k$-linear map $k[\overline{G}]\otimes A\rightarrow A$ (see \cite[\S 7.1]{Mont} for details). For $a\in A$ and $g\in \overline{G}$, we can identify them in $A\#_{\sigma}k[\overline{G}]$ with $a\#1$ and $1\#g$, respectively.

Let $V$ be any finite-dimensional subspace of $A\#_{\sigma}k[\overline{G}]$. Then it is easy to see that $V\subset A+Ag_1+\cdots+ Ag_m$ for some $g_1, \cdots, g_m\in \overline{G}$. Denote the set $\{g_1, g_2, \cdots, g_m\}$ by $S$ and let $g_S(n)$ be as defined in the previous lemma. Notice that $g_iA\subset Ag_i$ by $(\ref{product})$. As a consequence, 
$$V^n\subset\displaystyle (A+Ag_1+\cdots+ Ag_m)^{n}\subset  \sum _i Af_i,$$
where $f_i$ runs through all words of length $\le n$ in $\overline{G}$. Therefore, 
\begin{align*}
\limsup_{n\rightarrow \infty }\log_n\dim V^{n}&\le \limsup_{n\rightarrow \infty }\log_n(g_S(n)\cdot\dim A)\\
&=\limsup_{n\rightarrow \infty }\log_ng_S(n)\\
&\le \GK k[\overline{G}].
\end{align*}
Hence $\GK k[G]\le \GK k[\overline{G}]$, as desired.
\end{proof}

The proof of Lemma \ref{addition} is similar to the proof of \cite[Proposition 1]{LMO}. It seems that Lemma \ref{groupalgebra} is a corollary to Lemma \ref{addition}. However, there is some slight difference since we need to take the cocycle $\sigma$ into consideration. In fact, as we can see in the proof, the cocycle does not contribute to the GK-dimension because $N$ is finite. One might ask whether $\GK k[G]=\GK k[N]+\GK k[G/N]$ if $N$ is not finite. The answer is negative (see Example \ref{freenilpotent}).

\begin{lemma}\label{prime}
Let $G$ be a totally ordered group acting on a prime ring $A$. Then the skew group ring $B=A*G$ is also prime. 
\end{lemma}
\begin{proof}
For any non-zero elements $r, s\in B$, we need to show $rBs\ne 0$. Note that $B$ is $G$-graded. Let $ag$ and $cf$ be the leading terms of $r$ and $s$ respectively, where $a, c\in A$ and $g, f\in G$. Since $A$ is prime, there exists $b\in A$ such that $abc\ne 0$. Then 
$$(ag)(b^{g^{-1}}g^{-1})(cf)=abcf\ne 0.$$
Hence $r(b^{g^{-1}}g^{-1})s\ne 0$, which completes the proof.
\end{proof}

Now we are ready to prove the main theorem of this paper. We start with the case when $H$ is pointed and cosemisimple, i.e. $H\cong k[G]$ for some group $G$.
Before proving the proposition, we list some well-known facts about groups.

For a group $G$, we have the upper central series
$$\{1\}=Z_0(G)\subset Z_1(G)\subset Z_2(G)\cdots\subset Z_i(G)\subset\cdots,$$
where $Z_{i+1}(G)/Z_i(G)$ is the center of $G/Z_i(G)$. 
If $Z_{r-1}(G)\subsetneq Z_r(G)=G$ for some $r$, then $G$ is called a nilpotent group of nilpotency class $r$.

\begin{lemma}\label{nilpotent}
Let $G$ be a nilpotent group. Then the following statements are true.
\begin{enumerate}
\item[\textup{(1)}] Every subgroup of $G$ is nilpotent.
\item[\textup{(2)}] The elements of finite order form a normal subgroup $T$ such that $G/T$ is torsion-free.
\item[\textup{(3)}] If $G$ is finitely generated, then every subgroup of $G$ is also finitely generated.
\item[\textup{(4)}] The group $G$ is finite if it is generated by a finite number of elements each having finite order.
\item[\textup{(5)}] If $G$ is torsion-free, then $G/Z_i(G)$ is torsion-free for each $i\ge 0$.
\end{enumerate}
\end{lemma}
\begin{proof}
The first statement is easy. Statements $(2)$, $(3)$, $(4)$ and $(5)$ are \cite[Theorem 9.16]{Mac}, \cite[Corollary 9.18]{Mac}, \cite[Theorem 9.17]{Mac} and \cite[Lemma 8.2.3]{Ep} respectively.
\end{proof}

\begin{proposition}\label{cosemisimple}
Let $G$ be a group. If $2\le \GK k[G]< \infty$, then $G$ has a subgroup isomorphic to $\mathbb{Z}^2$.
\end{proposition}
\begin{proof}
By the definition of GK-dimension, 
$$\GK k[G]=\sup_B\GK B,$$
where $B$ runs through all finitely generated subalgebra of $k[G]$. It is easy to see that every finitely generated subalgebra $B$ of $k[G]$ is contained in some $k[L]$ where $L$ is a finitely generated subgroup of $G$. Hence 
$$\GK k[G]=\sup_L\GK k[L],$$
where $L$ runs through all finitely generated subgroups of $G$. By \cite[Theorem 11.1, 11.14]{KL}, $\GK k[L]$ must be an integer. So there is a finitely generated subgroup $N$ of $G$ such that $\GK k[N]=\GK k[G]$. Now by \cite[Theorem 11.1]{KL}, $N$ is nilpotent by finite, which means that $N$ has a nilpotent subgroup $P$ with $[N:P]<\infty$. Then by \cite[Proposition 5.5]{KL}, $\GK k[P]=\GK k[N]$. Also, by the same argument above, we can assume that $P$ is finitely generated. Therefore by replacing $G$ with $P$, we can assume that $G$ is finitely generated nilpotent.

First, we deal with the case that $G$ is torsion-free. If the nilpotency class $r=1$, i.e. $G$ is abelian, then $G$ must be a torsion-free abelian group of rank $\ge2$. Let $N$ be a subgroup of $G$ with rank $2$, then $N\cong\mathbb{Z}^2$. 
If $G$ has nilpotency class $r\ge 2$. Choose $h\in G\setminus Z_1(G)$. By Lemma \ref{nilpotent} part $(5)$, $\langle h\rangle \cap Z_1(G)=\{1\}$. Let $g\in Z_1(G)\setminus \{1\}$. Then by the choice of $h$, the subgroup $N$ generated by $g$ and $h$ is isomorphic to $\mathbb{Z}^2$.

For the general case, let $T$ be the set of torsion elements in $G$. By Lemma $\ref{nilpotent}$, $T$ is a finite normal subgroup of $G$ and $G/T$ is torsion-free. By Lemma \ref{groupalgebra}, we see that $\GK k[G/T] =\GK k[G]$. Now by the torsion-free case, there is a subgroup $B$ of $G$ such that $B/T\cong \mathbb{Z}^2$. By \cite[Lemma 8.2.4]{Ep}, $B$ has an abelian subgroup $N$ such that $[B: N]< \infty $. Then $\GK k[N]=\GK k[B]=\GK k[B/T]=2$, by \cite[Proposition 5.5]{KL} and Lemma $\ref{addition}$.  Hence $N$ is an abelian group of rank $2$, which means that it has a subgroup isomorphic to $\mathbb{Z}^2$. This completes the proof.
\end{proof}

Next we will deal with the case when $H$ is a connected Hopf algebra. Before that, we need a lemma on the existence of $2$-dimensional Lie subalgebra of a finite-dimensional Lie algebra. 
\begin{lemma}\label{2subalgebra}
Let $\mathbf{g}$ be a finite-dimensional Lie algebra. If $\dim_k \mathbf{g}\ge 2$, then $\mathbf{g}$ has a Lie subalgebra of dimension $2$.
\end{lemma}
\begin{proof}
If $\mathbf{g}$ is nilpotent, by \cite[3.2, Proposition]{H}, the center of $\mathbf{g}$ is non-zero. Choose a non-zero element $x$ in the center and another element $y$ which is not a scalar multiple of $x$. Then the subspace of $\mathbf{g}$ spanned by $x$ and $y$ is a $2$-dimensional Lie subalgebra.

If $\mathbf{g}$ is not nilpotent, by Engel's Theorem, there exists an element $x$ which is not ad-nilpotent. Hence the linear map $\ad x$ has a non-zero eigenvalue $\lambda$ with eigenvector $y$. Then the $k$-linear span of $x$ and $y$ is the desired Lie subalgebra.
\end{proof}

\begin{proposition}\label{connected}
Let $H$ be a connected Hopf algebra. If $2\le \GK H< \infty$, then $H$ contains a Hopf subalgebra isomorphic to $U(\mathbf{h})$, where $\mathbf{h}$ is a $2$-dimensional Lie algebra.
\end{proposition}
\begin{proof}
Let $K$ be the Hopf subalgebra generated by $P(H)$, the set of primitive elements in $H$. By \cite[Theorem 5.6.5]{Mont}, $K\cong U(\mathbf{g})$, where $\mathbf{g}=P(H)$. First we claim that $\dim_k \mathbf{g}\ge 2$.

If $\dim_k \mathbf{g}=1$, then $K\cong k[y]$, the polynomial ring of one-variable, as Hopf algebras. Obviously, $K\subsetneq H$ since they have different GK-dimensions. Let $m$ be the smallest number such that $K_{m}\subsetneq H_{m}$. By the construction of $K$, we have $m\ge 2$. Choose $z\in H_{m}\backslash K_{m}$. By \cite[Lemma 5.3.2]{Mont} and its proof, we can assume that
$$\Delta(z)=1\otimes z+ z\otimes 1+ w,$$
where $w\in \sum\limits_{i=1}^{m-1}H_{m-i}\otimes H_i=\sum\limits_{i=1}^{m-1}K_{m-i}\otimes K_i$. Moreover, by replacing $z$ with $z-\e(z)\cdot 1$, we can assume $\e(z)=0$. Notice that $K\cong k[y]$ is coradically graded with the usual grading and therefore $K_n$ is spanned by all polynomials in $y$ of degree $\le n$.  Also, the grading on $K$ induces a grading on $K^{\otimes \ell }$ by setting $\deg f_1\otimes f_2\otimes \cdots \otimes f_{\ell }=\deg f_1+\cdots+ \deg f_{\ell }$. Hence we can write
$$w=\sum\limits_{i+j\le m }\alpha_{ij}y^i\otimes y^j.$$
By the counit axiom,
$$z=(Id\otimes \e)\Delta(z)=z+\sum \limits_{i\le m }\alpha_{i0}y^i,$$
$$z=(\e\otimes Id)\Delta(z)=z+\sum \limits_{j\le m }\alpha_{0j}y^j. $$
Hence $\alpha_{i0}=0$ and $\alpha_{0j}=0$ for any $i, j\le m$.

If $m=2$, then
$$\Delta(z)=1\otimes z+ z\otimes 1+ \alpha_{11}y\otimes y,$$
Let $z'=z-\frac{\alpha_{11}}{2}y^2$. Then $z'\in H_{2}\backslash K_{2}$ by the choice of $z$. But a direct calculation shows that
$$\Delta(z')=1\otimes z'+ z'\otimes 1.$$
So $z'\in H_{1} =K_{1}$, which is a contradiction.

Now assume that $m\ge 3$. Define two maps
$$\partial ^1: H\rightarrow H\otimes H\,\,\,\,\,\text{and}\,\,\,\,\,\partial^2: H\otimes H\rightarrow H\otimes H\otimes H$$
by $\partial^1(c)=1\otimes c-\Delta (c)+c\otimes 1$ and $\partial^2(c\otimes d)=1\otimes c\otimes d- \Delta(c)\otimes d+c\otimes \Delta(d)-c\otimes d\otimes 1$ for any $c, d\in H$. By \cite[Corollary XVIII.5.2]{Ka}, $\partial^2\circ\partial ^1=0$. Hence $\partial^2(w)=0$. 

Let $u=\sum\limits_{i+j=m }\alpha_{ij}y^i\otimes y^j$.
Notice that the map $\partial^2$, when restricted from $K\otimes K$ to $K\otimes K\otimes K$, preserves the grading. Therefore $\partial^2(u)=0$. Also, by replacing $z$ with $z-\beta y^m$, where $\beta=\alpha_{m-1, 1}/m$, we can assume that $\alpha_{m-1,1}=0$. Now, 
\begin{align*}
\partial^2(u)=&\partial ^2(\sum \limits_{i=2}^{m-1}\alpha_{m-i,i} y^{m-i}\otimes y^i)\\
=&\nonumber 1\otimes \sum \limits_{i=2}^{m-1}\alpha_{m-i,i}  y^{m-i}\otimes y^i-\sum \limits_{i=2}^{m-1}\alpha_{m-i,i}  \Delta(y^{m-i})\otimes y^i\\
&+\sum \limits_{i=2}^{m-1}\alpha_{m-i,i}  y^{m-i}\otimes \Delta(y^i)-\sum \limits_{i=2}^{m-1}\alpha_{m-i,i}  y^{m-i}\otimes y^i\otimes 1\nonumber\\
=&-\sum \limits_{i=2}^{m-1}  \sum \limits_{j=0}^{m-i-1}\alpha_{m-i,i}  \binom{m-i }{j}y^{m-i-j}\otimes y^j\otimes y^i
+\sum \limits_{i=2}^{m-1}\sum \limits_{j=1}^{i}\alpha_{m-i,i} \binom{i}{j}y^{m-i}\otimes y^{i-j}\otimes y^j.\nonumber\\
\end{align*}
Since $\partial^2(u)=0$, by comparing the coefficients of the terms $y^{m-i}\otimes y^{i-1}\otimes y$, we see that $\alpha_{m-i,i} =0$ for $i=2, \cdots, m-1$. It then follows that 
$$\Delta(z)=1\otimes z+z\otimes 1+\sum\limits_{i,j\ge 1\atop i+j\le m-1 }\alpha_{ij}y^i\otimes y^j.$$
But this implies that $\Delta(z)\in H_0\otimes H+ H\otimes H_{m-2}$, i.e. $z\in H_{m-1}=K_{m-1}$, which is a contradiction.

Now $K\cong U(\mathbf{g})$ where $\mathbf{g}$ is a finite-dimensional Lie algebra and $\dim_k \mathbf{g}\ge 2$. By Lemma \ref{2subalgebra}, $\mathbf{g}$ has a $2$-dimensional Lie subalgebra $\mathbf{h}$. Then $U(\mathbf{h})$ is a Hopf subalgebra of $K\cong U(\mathbf{g})$, and thus a Hopf subalgebra of $H$. This completes the proof.
\end{proof}

Finally, we have to handle the case when $H$ is neither cosemisimple nor connected.
\begin{lemma}\label{case3}
Let $H$ be a pointed Hopf algebra with finite GK-dimension. If $H$ is a domain and is not connected, not cosemisimple, then $H$ contains a Hopf subalgebra isomorphic to either $A(b, \xi)$ or $C(b+1)$ for some integer $b$ and $\xi \in k^{\times }$.
\end{lemma}
\begin{proof}
Since $H$ is not cosemisimple, by \cite[5.4.1]{Mont}, we can find $g, y\in H$ such that $g$ is group-like and $y\notin H_0$ is $(1, g^b)$-primitive for some integer $b$. Notice that $g$ is torsion-free since $H$ is a domain. Moreover, we can choose $g\neq 1$ since $H$ is not connected. By Corollary $\ref{existence2}$,  there are $b_0, \cdots, b_m, \beta\in k$ and $\xi \in k^{\times }$ such that $f=b_0y+b_1gyg^{-1}+\cdots+b_{m}g^{m}yg^{-m}\in  H_1\backslash H_0$ and $gf=\xi fg+\beta g(g^b-1)$.

Let $K$ be the Hopf subalgebra generated by $f$ and $g^{\pm 1 }$. By construction, $K$ is neither cosemisimple nor connected.
If $\xi\ne1$,  by replacing $f$ with $f'=f-\frac{\beta }{1-\xi }(g^b-1)$, we have $gf'=\xi f'g$. Hence $K$ is a Hopf quotient of $A(b, \xi)$. If $\xi =1$ and $\beta\ne 0$,  by replacing $f$ with $f'=\beta^{-1}f$, we have $gf'=f'g+ g(g^b-1)$, which means that $K$ is a Hopf quotient of $C(b+1)$. If $\beta=0$, $K$ is a Hopf quotient of $A(b, \xi)$. In all the above cases, we can apply Lemma $\ref{AandC}$, which yields that $K\cong A(b, \xi)$ or $K\cong C(b+1)$.
\end{proof}

Combining the previous lemmas together, we have the main theorem of this paper.

\begin{theorem}\label{GK2sub}
Let $H$ be a pointed Hopf algebra over an algebraically closed field $k$ of characteristic $0$. If $2\le \GK H<\infty$ and $H$ is a domain, then $H$ contains a Hopf subalgebra of GK-dimension $2$. To be precise, 
\begin{enumerate}
\item[\textup{(1)}]  If $\GK k[G]=0$, where $G=G(H)$, then $H$ contains a Hopf subalgebra isomorphic to $U(\mathbf{h})$, where $\mathbf{h}$ is a $2$-dimensional Lie algebra.
\item[\textup{(2)}]  If $\GK k[G]\ge 2$, then $H$ contains a Hopf subalgebra isomorphic to $k\mathbb{Z}^2$.
\item[\textup{(3)}]  If $H_1\supsetneq H_0=k[G]$ and $\GK k[G]\ge 1$, then $H$ contains a Hopf subalgebra isomorphic to either $A(b, \xi)$ or $C(b+1)$ for some integer $b$ and $\xi \in k^{\times }$.
\end{enumerate}
\end{theorem}
\begin{proof}
Since $\GK H<\infty$, we have $\GK k[G]<\infty$.

If $\GK k[G]=0$, then by assumption, $G$ is trivial. Hence $H$ is connected. By Proposition $\ref{connected}$, $H$ contains a Hopf subalgebra isomorphic to $U(\mathbf{h})$, where $\mathbf{h}$ is a $2$-dimensional Lie algebra.

If $\GK k[G]\ge 2$, then by Proposition $\ref{cosemisimple}$, $k[G]$, and thus $H$, contains a Hopf subalgebra isomorphic to $k\mathbb{Z}^2$.

Part (3) follows from Lemma $\ref{case3}$.
\end{proof}

\begin{remark} 
It seems that in Theorem \ref{GK2sub}, if the condition that $H$ is a domain is dropped, we can still find Hopf subalgebra of GK-dimension 2.
By Proposition $\ref{cosemisimple}$ and Proposition $\ref{connected}$, this is true if $H$ is either a group algebra or connected. However, we still do not know how to deal with the case when $H$ is neither a group algebra nor connected.
\end{remark}

One might ask whether Lemma $\ref{case3}$ still holds if $H$ is a prime ring instead of a domain. The answer is negative as shown in the following example. 

\begin{example}\label{counter}
Let $E=k\langle x_0^{\pm 1}, \cdots, x_n^{\pm 1 }, y\mid yx_i+x_iy=0, y^2=x_0^2-1, x_ix_j-x_jx_i=0\rangle$. Then $E$ has a Hopf algebra structure under which $x_i$ is group-like for $i=0, 1, \cdots, n$ and $y$ is $(1, x_0)$-primitive. Denote this Hopf algebra by $E(n)$. It is easy to check that $E(n)$ is not connected, not cosemisimple. Notice that $E(n)$ has a basis consisting of elements of the form $x_0^{i_0}\cdots x_n^{i_n}y^\epsilon$, where $\epsilon$ is either $0$ or $1$. Since $y^\epsilon\in E(n)_1$, the element $x_0^{i_0}\cdots x_n^{i_n}y^\epsilon$ is in $E(n)_1$ as well. Hence, $E(n)=E(n)_1$. By \cite[Lemma 5.2.12]{Mont}, every Hopf subalgebra $H$ of $E(n)$ has $H=H_1$. However, for $D=A(b, \xi)$ or $C(b+1)$, it is easy to find elements in $D_2\backslash D_1$. So $E(n)$ has no Hopf subalgebra isomorphic to either $A(b, \xi)$ or $C(b+1)$.

If $n=0$, the statement follows from \cite[Lemma 2.3 and 2.6]{G}.
If $n\ge 1$, then there is a Hopf projection from $E(n)$ to $H=k[x_1^{\pm1 }, \cdots, x_n^{\pm 1 }]\cong k[\mathbb{Z}^n]$. Hence by \cite[Proposition 7.2.3]{Mont}, $E(n)\cong A\#H$, where $A=E(n)^{co H }$. In fact, $A$ is the subalgebra of $E(n)$ generated by $y$ and $x_0^{\pm 1 }$. So $A$ is a Hopf subalgebra isomorphic to $E(0)$, which is prime of GK-dimension $1$. Then it follows from Lemma \ref{addition} and \ref{prime} that $E(n)$ is prime with $\GK E(n)=n+1$.
\end{example}

The next example shows that Theorem \ref{GK2sub} fails if the Hopf algebra $H$ is not pointed.
\begin{example}
Let $H$ be $\mathcal{O}(SL_2(k))$, the coordinate ring of the algebraic group $SL_2(k)$. It is well known that $H$ has a $4$-dimensional simple subcoalgebra, so $H$ is not pointed. By \cite[Theorem I.2.10]{BG}, $H$ is a noetherian domain. Also, $\GK H=3$ by \cite[Corollary]{LS}.

Now, $H$ is a commutative noetherian domain of GK-dimension $3$. Moreover, by \cite[\S 2.4]{Mont}, $H$ is cosemisimple. So any subcoalgebra of $H$ is also cosemisimple by \cite[Lemma 5.2.12]{Mont}. If $H$ has a Hopf subalgebra $K$ of GK-dimension $2$, then by \cite[Proposition 2.3]{GZ},  $K$ must be isomorphic to an enveloping algebra $U(\mathbf{g})$, where $\mathbf{g}$ is 2-dimensional abelian, or a group algebra $k\Gamma$, where $\Gamma$ is free abelian of rank $2$, or a Hopf algebra $A(b, 1)$ for some nonnegative integer $b$. Since $U(\mathbf{g})$ and $A(b, 1)$ are not cosemisimple, $K$ is of the form $k\Gamma$. But the only invertible elements in $H$ are just non-zero scalars. Therefore, $H$ has no Hopf subalgebras of GK-dimension $2$.
\end{example}
Another question is whether a Hopf algebra with GK-dimension $>3$ contains a Hopf subalgebra of GK-dimension three. The answer is negative. In fact, this statement is not even true for group algebras, as shown in the following example.

\begin{example}\label{freenilpotent}
Let $G$ be the free nilpotent group of class $2$ on generators $x$ and $y$. Thus, if we write $z: =[x, y]$, then 
$$[x, z]=[y,z]=1.$$
It is easy to check that $G_1=Z(G)=\langle z\rangle$. By \cite[Theorem 11.14]{KL}, $\GK k[G]=4$. 
By construction, $G/G_1\cong \mathbb{Z}^2$. Hence
$$\GK k[G_1]+\GK k[G/G_1]=3<4=\GK k[G],$$
which gives a negative answer to the question given after Lemma $\ref{groupalgebra}$.

Notice that every Hopf subalgebra of $k[G]$ must be of the form $k[L]$ where $L$ is a subgroup of $G$. Therefore, to show $k[G]$ has no Hopf subalgebra of GK-dimension three, we only need to show that $G$ has no subgroup of cubic growth. Suppose to the contrary that $N$ is a subgroup of $G$ with cubic growth.

If $N$ has nilpotency class $1$, then $N$ must be a finitely generated torsion-free abelian group of rank $3$, i.e. $N\cong \mathbb{Z}^3$. Notice that we have an exact sequence of groups,
$$1\rightarrow N\cap G_1\rightarrow N\rightarrow N/N\cap G_1\rightarrow 1.$$
If $N\cap G_1=\{1\}$, then $N$ embeds in $G/G_1$, which is impossible since $G/G_1$ has only rank $2$. Hence $N\cap G_1$ has rank $1$, which implies that $N/N\cap G_1$ has rank $2$ and $G_1/N\cap G_1\cong (N+G_1)/N$ is finite. By considering the exact sequence
$$1\rightarrow N/N\cap G_1\rightarrow G/G_1\rightarrow  G/(N+G_1)\rightarrow 1,$$
we see that $G/(N+G_1)$ is also finite. It follows that $G/N$ is finite. Then by \cite[Proposition 5.5]{KL}, $\GK k[G]=\GK k[N]$, which is a contradiction.

Now suppose $N$ has nilpotency class $2$. Then $N_1$ is a non-trivial subgroup of $G_1$, so it must have rank $1$. By \cite[Theorem 11.14]{KL}, $N/N_1$ has rank $1$. Choose $w\in N\backslash N_1$ such that the image of $w$ in $N/N_1$ is torsion-free. Denote $\langle w, N_1\rangle$ by $P$, then $N/P$ is finite. On the other hand, $w$ commutes with $N_1$. Hence by construction, $P$ is abelian of rank $2$. Now $\GK N=\GK P=2$ since $N/P$ is finite. But this again is a contradiction.
\end{example}

\section{Acknowledgements}
The author thanks Professor James Zhang, Professor Ken Brown, Xingting Wang and Cris Negron for useful conversations and for their careful reading of this paper. The author would also like to thank the referee for the valuable comments.

\bibliographystyle{model1-num-names}
\bibliography{<your-bib-database>}

\end{document}